\documentclass[11pt]{article}

\usepackage[latin1]{inputenc}

\usepackage{amsfonts,amsmath,amssymb,amsthm,graphicx,epsfig,float}
\usepackage[T1]{fontenc}

\usepackage[english,francais]{babel}

{
  \newtheorem{theoreme}{Th\'eor\`eme}
  \newtheorem*{theoreme*}{Th\'eor\`eme}
  \newtheorem{lemme}[theoreme]{Lemme}

  \newtheorem{proposition}[theoreme]{Proposition}
\newtheorem*{corollaire*}{Corollaire}
\newtheorem*{proposition*}{Proposition}
\theoremstyle{remark}
  \newtheorem*{remarque*}{Remarque}
}

\newcounter{ex}

\newenvironment{rem*}{
  \noindent\textbf{Remarque. }}{}

\newcommand{\Cc}{\mathbb{C}}
\newcommand{\Nn}{\mathbb{N}}

\newcommand{\Pp}{\mathbb{P}}

\renewcommand {\epsilon}{\varepsilon}

\renewcommand {\leq}{\leqslant}
\renewcommand {\geq}{\geqslant}

\title{{\bf Un phénomène de concentration de genre}}
\author{Henry de Thélin}
\date{}

\begin{document}
\maketitle


\def\figurename{{Fig.}}%
\def\proofname{Preuve}
\def\contentsname{Sommaire}%

\begin{abstract}
Nous montrons que le genre des préimages d'une droite projective par un endomorphisme holomorphe générique de $\Pp^2(\Cc)$ se concentre sur le support de la mesure de Green.
\end{abstract}

\selectlanguage{english}
\begin{center}
{\bf{A phenomenon of concentration of genus }}
\end{center}

\begin{abstract}
We show that the main part of the genus of the preimages of a projective line by a generic holomorphic endomorphism of $\Cc \Pp^2$ goes to the support of the Green measure.
\end{abstract}

\selectlanguage{francais}

Mots-clefs: dynamique holomorphe, laminarité.\\
AMS: 32H50, 32U40, 32Q45.

\section*{{\bf Introduction}}
\par
E. Bedford, M. Lyubich et J. Smillie ont introduit dans \cite{BLS} une notion de lamination en un sens faible qui décrit les ensembles de Julia des applications de Hénon: les courants laminaires. Ce sont des $(1,1)$-courants positifs qui s'écrivent localement comme une intégrale de courants d'intégration sur une famille de disques holomorphes disjoints, hors d'un ensemble négligeable.\\
Les courants laminaires peuvent s'obtenir comme limite de courbes
analytiques dont le genre n'augmente pas plus vite que l'aire (voir
\cite{BLS}, \cite{Du1} et \cite{Det}). Plus précisément, si la
situation est locale on a:
\begin{theoreme*}{\label{laminarite}}{Voir \cite{Det}.}\\
Soit $C_n$ une suite de courbes analytiques lisses de la boule unité $B$ de $\Cc^2$.\\
On note $A_n$ l'aire de $C_n$, $G_n$ le genre de $C_n$ et on suppose que $T_n=\frac{[C_n]}{A_n}$ converge vers un $(1,1)$-courant positif fermé $T$ de $B$.\\
Alors, si $G_n=O(A_n)$, $T$ est laminaire.
\end{theoreme*}
On peut utiliser ce théorème pour l'étude de courants limites de $\frac{[C_n]}{A_n}$ où $C_n$ est une courbe algébrique lisse de $\Pp^2(\Cc)$ (par exemple $C_n=f^{-n}(L)$ où $f$ est un endomorphisme holomorphe de $\Pp^2(\Cc)$ et $L$ une droite projective générique). En effet, malgré un genre total en $O(A_n^2)$, si on sait trouver un ouvert où le genre se concentre peu (en $O(A_n)$), on en déduit la laminarité des courants limites dans celui-ci. Cela nous conduit naturellement à l'objet de cet article: la concentration du genre des courbes lisses précédentes, issues de la dynamique.
\par
A partir d'un endomorphisme holomorphe, $f$, de degré $d \geq 2$,
J.E. Forn{\ae}ss et N. Sibony ont défini le courant de Green, $T$,
associé à $f$ (voir \cite{FS}), dont le support est l'ensemble de
Julia de $f$. Ce courant possède un potentiel continu: on peut donc
définir son auto-intersection $\mu=T \wedge T$ (voir
\cite{FS}). D'autre part, en généralisant un résultat de Forn{\ae}ss
et Sibony (voir \cite{FS2}), C. Favre et M. Jonsson (voir \cite{FJ})
ont montré que ce courant est naturel d'un point de vue dynamique: il
équirépartit les préimages de droites génériques. En effet, si
$C_n=f^{-n}(L)$ avec $L$ une droite projective, alors
$\frac{[C_n]}{A_n}=\frac{[C_n]}{d^n}$ converge vers $T$ sauf pour un
ensemble algébrique de droites $L$. Autrement dit, dans notre situation, l'étude de la laminarité de $T$ se ramène à celle du genre de $f^{-n}(L)$.\\\\
Quand $f$ est un endomorphisme critiquement fini, i.e. dont le lieu critique est prépériodique, cette approche aboutit:
\begin{theoreme}{\label{theo1}}
Pour un endomorphisme holomorphe $f$ critiquement fini, le genre de
l'image réciproque par $f^n$ d'une droite projective générique $L$,
$f^{-n}(L)$, est dominé par $O(d^n)$ en dehors d'un voisinage du
support de la mesure $\mu$. Ici $L$ générique signifie que $L$ est en
dehors d'une union dénombrable d'ensembles algébriques du dual de $\Pp^2(\Cc)$.

\end{theoreme}

\begin{corollaire*}

Pour un endomorphisme holomorphe critiquement fini, le courant de Green $T$ est laminaire en dehors du support de la mesure $\mu$.

\end{corollaire*}

Remarquons qu'il existe des exemples d'endomorphismes critiquement
finis dont le courant de Green est laminaire seulement en dehors du
support de $\mu$. Il suffit de prendre la suspension à $\Pp^2(\Cc)$
d'un exemple de Lattès de $\Pp^1(\Cc)$. Par ailleurs, il est naturel
de ne pas avoir de la laminarité sur le support de $\mu$. En effet dans
\cite{Du2}, R. Dujardin a montré qu'un courant $T$ fortement
approximable (condition un peu plus forte que laminaire) sur le
support de $T \wedge T$, à potentiel continu vérifie $T \wedge
T=0$. Signalons enfin qu'il existe des
endomorphismes critiquement finis pour lesquels le support de $\mu$
est différent du support de $T$ (par exemple $f([z:w:t])=[z^d:w^d:t^d]$).\\

Pour $f$ générique, le genre des préimages d'une droite se concentre toujours sur le support de $\mu$: hors de celui-ci, on obtient un contrôle en $O(d^{n(1+\epsilon)})$. Plus précisément, le résultat principal de ce texte est le:

\begin{theoreme}{\label{genre}}
Pour $f$ générique parmi les endomorphismes de degré $d$, on a:
$$\limsup_{n \rightarrow \infty} \frac{1}{n} \log \max_{L \in
  (\Pp^2)^{*}} \mathrm{Genre }(f^{-n}(L) - U )\leq \log d,$$
où $U$ est un petit voisinage du support de $\mu$.
\end{theoreme}

La démonstration de ce théorème se fera essentiellement en deux étapes. En effet, si on admet un instant que les anses de $f^{-n}(L)$ sont infiniment petites, on constate qu'une anse de $f^{-1}(L)$, tirée en arrière par $f^{n-1}$, se comporte comme $f^{-(n-1)}(x)$ (où $x$ est un point de $\Pp^2(\Cc)$). Autrement dit, la première étape de la démonstration consistera à contrôler le nombre de points de $f^{-n}(x)$ hors d'un petit voisinage du support de $\mu$. Celle-ci s'énonce:
\begin{proposition}{\label{prop1}}
Pour $f$ générique parmi les endomorphismes de degré $d$, on a:
$$\limsup_{n \rightarrow \infty} \frac{1}{n} \log \max_{x \in \Pp^2}
\mathrm{Cardinal }(f^{-n}(x) - U) \leq \log d,$$
où $U$ est un petit voisinage du support de $\mu$.

\end{proposition}

La seconde étape consistera alors à dominer la taille des anses. Pour cela, on utilisera des modules d'anneaux et des comparaisons aire-longueur.\\\\
Voici donc le plan de ce texte: dans le premier paragraphe, on traitera le cas critiquement fini tandis que dans le second on démontrera la proposition \ref{prop1}. Enfin la dernière partie sera consacrée à la preuve du théorème \ref{genre}.\\\\
{\bf Remerciement:} Je tiens à remercier mon directeur de thèse Julien Duval pour son aide précieuse dans l'élaboration de cet article.

\section{{\bf Le cas critiquement fini}}{\label{par1}}

Un endomorphisme holomorphe est critiquement fini si son lieu postcritique, $\mathcal{C}=\cup_{n \geq 0} f^n(C_f)$ (où $C_f$ est l'ensemble critique de $f$), est une courbe algébrique (voir \cite{FS1}).\\

{\textit{Démonstration du théorème \ref{theo1}.}}\\
Pour calculer le genre de $f^{-n}(L)$, on va construire une
triangulation de $L$ qui se relèvera en une triangulation de
$f^{-n}(L)$ (par triangulation, on entend décomposition en des
disques). Pour conclure, il restera alors à faire un calcul de
caractéristique d'Euler sur $f^{-n}(L)$ et à utiliser le fait que le
genre d'une surface connexe $M$ qui a $b$ composantes de bord est égal à $1- \frac{\chi(M)+b}{2}$.\\\\
On va supposer que $f^{-n}(L)$ est lisse (ce qui est vrai pour $L$ générique). D'autre part, on va traiter le cas d'une droite $L$ qui passe par un point $a$ appartenant au support de $\mu$ privé de $\mathcal{C}$ telle que $\mbox{Cardinal}( L \cap \mathcal{C})=\mbox{degré de } \mathcal{C}=\tau$ (c'est possible car $\mu$ ne charge pas $\mathcal{C}$). Le cas d'une droite générique se fait d'une manière proche.\\
Dans la suite, on supposera que $a$ est le pôle nord de $L$.\\
Quitte à bouger un peu la droite, on peut supposer que l'intersection des méridiens qui passent par $\mathcal{C} \cap L$ avec l'équateur est constituée de $\tau$ points: $a_1,...,a_{\tau}$. On les ordonne via une orientation de l'équateur.\\
En considérant alors les méridiens qui passent par le milieu des segments $[a_i,a_{i+1}]$ ($i=1,...,\tau -1$) et $[a_{\tau},a_1]$ (où les segments considérés sont situés sur l'équateur), on obtient une triangulation de $L$ avec $\tau$ faces, $\tau$ arêtes et deux sommets.\\
Cette triangulation se relève par $f^n$ en une triangulation de $f^{-n}(L)$. En effet, un disque qui rencontre $\mathcal{C}$ en au plus un point se relève en un disque par $f^n$.\\
On va maintenant calculer le genre de $f^{-n}(L)$ hors d'un $\gamma$-voisinage du support de la mesure $\mu$, $U_{\gamma}$.\\ 
Soit $\Gamma$ la réunion des arêtes et des sommets de la triangulation qui sont contenus dans $U_{\gamma}$;\\
alors $g(f^{-n}(L)-U_{\gamma}) \leq g(f^{-n}(L)-\Gamma)$.\\
D'autre part, une composante connexe de $f^{-n}(L)- \Gamma$ peut être de deux types:\\
- soit elle ne contient pas d'arête de la triangulation et alors c'est un disque (donc de genre nul),\\
- soit elle en contient au moins une.\\
Si $\mathcal{Q}$ désigne l'union des composantes du deuxième type, on obtient :
$$g(f^{-n}(L) - U_{\gamma}) \leq \mbox{Nombre de composantes de } \mathcal{Q} - \frac{ \chi(\mathcal{Q})}{2}.$$
Ainsi:
$$ g(f^{-n}(L) - U_{\gamma}) \leq \frac{3}{2} \mbox{du nombre d'arêtes qui sortent de }U_{\gamma}.$$
En effet les faces qui composent les éléments de $\mathcal{Q}$ étant des disques, $- \chi(\mathcal{Q})$ est majoré par le nombre d'arêtes qui sortent de $U_{\gamma}$.\\
Il reste donc à majorer ce dernier terme.\\

Les arêtes étant attachées au support de $\mu$, on en déduit que le genre de $f^{-n}(L)$ hors de $U_{\gamma}$ est majoré (à une constante près) par le nombre d'arêtes de la triangulation de diamètre supérieur à $\gamma$. Il suffit donc de majorer ce terme par $O(d^n)$. C'est l'objet du:

\begin{lemme} 
Les préimages par $f^n$ des arêtes de la triangulation de $L$ sont de diamètre inférieur à $\gamma$ (sauf $O(d^n)$ d'entre elles).

\end{lemme}

\begin{proof}

La démonstration de ce lemme repose sur la comparaison aire-diamètre qui suit (voir \cite{BD2}).\\
\underline{Fait}:\\
Il existe $C > 0$ tel que, pour toute paire de disques holomorphes $D \subset \tilde{D}$ dans $\Pp^2(\Cc)$, on ait
$$(\mbox{Diam}(D))^2 \leq C \frac{\mbox{Aire}(\tilde{D})}{\min(1,\mbox{Mod}(A))},$$
où $A$ désigne l'anneau $\tilde{D}-D$.\\

On peut mettre une arête $\alpha$ de la triangulation de $L$ dans un disque $D$ disjoint de $\mathcal{C}$, de sorte que le module de l'anneau $D - \alpha$ soit égal à une constante $m$ inférieure à $1$ et qui ne dépend que de $\mathcal{C} \cap L$.\\
On note $f_i^{-n}$ les branches inverses de $f^n$ sur $D$. Les anneaux $f_i^{-n}(D - \alpha)$ ont le même module que celui de $D - \alpha$.\\
En utilisant le fait précédent, on obtient donc une majoration du carré du diamètre de $f_i^{-n}(\alpha)$ par $ \frac{C}{m} \mbox{Aire}(f_i^{-n}(D))$. Autrement dit, le nombre de branches inverses pour lesquelles le diamètre de $f_i^{-n}(\alpha)$ est supérieur à $\gamma$ est majoré par $\frac{C}{m \gamma^2} \mbox{Aire}(f^{-n}(L)) = \frac{Cd^n}{m \gamma^2}$. C'est ce que l'on voulait démontrer.

\end{proof}
\qed \\
Dans le cas d'un endomorphisme quelconque, si on découpe $L$ avec des disques qui ne rencontrent les valeurs critiques $V=f(C_f)$ qu'en un seul point et que l'on remonte cette triangulation par $f$, on en obtient une de $f^{-1}(L)$. Cependant les éléments de cette triangulation n'ont aucune raison de continuer à rester des disques quand on les relèvera de nouveau par $f$ (certains d'entre eux peuvent toucher $V$ en plusieurs points). Pour conserver une triangulation, il faudra donc faire un redécoupage.\\
Lors de cette opération, on pourra garder le contrôle des longueurs des préimages des arêtes (on aura au plus $d^{n(1+ \epsilon)}$ arêtes de longueur supérieure à $\gamma$). On compensera alors la perte de l'attache par le contrôle en $d^{n(1+ \epsilon)}$ du nombre de points de  $f^{-n}(x)$ hors de $U_{\frac{\gamma}{2}}$, pour tout point $x$ de $\Pp^2(\Cc)$ (ce contrôle sera vrai pour des endomorphismes $f$ génériques). Cela suffira pour majorer le genre car les arêtes qui sortent de $U_{\gamma}$  sont composées des arêtes de longueur supérieure à $\gamma$, auxquelles on ajoute celles qui ont un sommet hors de $U_{\frac{\gamma}{2}}$.\\\\
Le plan de la démonstration du théorème \ref{genre} sera donc le suivant: dans un premier temps on montrera que l'on peut contrôler le nombre de points de $f^{-n}(x)$ hors de $U_{\frac{\gamma}{2}}$ par $d^{n(1+ \epsilon)}$. Puis, on verra que cela permet effectivement de majorer le genre de $f^{-n}(L)$ hors de $U_{\gamma}$ par $d^{n(1+ \epsilon)}$ .

\section{{\bf Contrôle des préimages des points}}

Le contrôle du nombre d'antécédents éloignés du support de $\mu$ d'un point $x$ de $\Pp^2(\Cc)$ s'inspire du calcul de l'entropie topologique d'un endomorphisme holomorphe de degré $d$.\\
Après quelques rappels sur l'entropie topologique qui incluront les différentes étapes de ce calcul, on passera au contrôle des préimages de points.

\subsection{{\bf Entropie topologique}}{\label{entropie}}

L'entropie topologique $h_{\mbox{top}}(f)$ est définie par (voir par exemple \cite{KH}):
$$h_{\mbox{top}}(f)=\sup_{\delta > 0} \limsup_{n \rightarrow \infty} \frac{1}{n} \log (\max \{\mbox{Card}(F),\mbox{ } F \mbox{ } (n,\delta) \mbox{-séparé} \})$$
où un ensemble est dit $(n,\delta) \mbox{-séparé}$ si pour tout couple $(x,y) \in F^2$ on a $d_n(x,y):=\max_{0 \leq q \leq n-1} d(f^q(x),f^q(y))  \geq \delta$.\\
C'est donc une quantité qui décrit le nombre d'orbites que l'on peut discerner à une erreur $\delta$ près en un temps $n$.\\
Un endomorphisme $f$ holomorphe de $\Pp^2(\Cc)$ est d'entropie topologique $2 \log d$. L'obtention de cette valeur est la combinaison d'une  majoration obtenue par M. Gromov (voir \cite{G}) et d'une minoration due à M. Misiurewicz et F. Przytycki (voir \cite{KH}).\\

La majoration repose de manière cruciale sur le théorème de Lelong (voir \cite{L}).\\
Si on note $\Gamma_n=\{ (x,f(x),...,f^{n-1}(x)), x \in \Pp^2(\Cc) \}$ le multigraphe de $f$ d'ordre $n$, on voit qu'un ensemble $F$ $(n,\delta) \mbox{-séparé}$ donne un ensemble $G$ $\delta$-séparé dans $\Gamma_n$ pour la distance produit (qui est $d_n$). En désignant par $\omega_n$ la forme kählérienne sur $(\Pp^2(\Cc))^n$ induite par la forme de Fubini-Study $\omega$ sur chaque facteur, on a:
$$\int_{\Gamma_n} \omega_n^2 = \mbox{vol}(\Gamma_n) \geq \sum_{y \in G} \mbox{vol}(B_n(y, \frac{\delta}{2}) \cap \Gamma_n )$$
où $B_n(y, \frac{\delta}{2})$ est la boule centrée en $y$ de rayon $\frac{\delta}{2}$ pour la métrique $d_n$. Le théorème de Lelong nous conduit alors à une minoration du volume de $\Gamma_n$ par $C (\mbox{Card}(G))$ (où $C$ est une constante indépendante de $n$). Autrement dit, l'entropie topologique est majorée par la quantité:
$$lov(f):= \limsup_{n \rightarrow \infty} \frac{1}{n} \log(\mbox{vol}(\Gamma_n)) = \limsup_{n \rightarrow \infty} \frac{1}{n} \log \int_{\Gamma_n} \omega_n^2$$
qui est égale à 
$$\limsup_{n \rightarrow \infty} \frac{1}{n} \log \int_{\Pp^2(\Cc)} \sum_{i,j=0}^{n-1} f^{i *} \omega \wedge f^{j *} \omega  = 2 \log d,$$
car $f^{i *} \omega$ est cohomologue à $d^i \omega$.\\\\

La minoration repose sur la construction d'ensembles $(n,\delta) \mbox{-séparés}$ dans les préimages $f^{-n}(x)$ de points $x$ dont les antécédents ne s'approchent pas trop souvent de l'ensemble critique.\\
Elle s'inscrit dans un cadre plus général, celui des applications de classe $C^1$:

\begin{theoreme}{(Misiurewicz-Przytycki)}\\
Si $M$ est une variété lisse, compacte et orientable et $f: M \rightarrow M$ une application de classe $C^1$, on a:
$$h_{\mathrm{top}}(f) \geq \log | \mathrm{deg}(f)|.$$

\end{theoreme}

En voici le schéma (voir \cite{KH}).\\
Supposons tout d'abord que $f$ est un revêtement de degré $N$.\\
Pour $\delta$ assez petit, tous les points de $f^{-1}(x)$ sont à distance au moins $\delta$ les uns des autres (et ceci uniformément en $x$). De là, on en déduit que les points de $f^{-n}(x)$ forment un ensemble $(n,\delta) \mbox{-séparé}$ de cardinal $N^n$ (pour $x$ fixé dans $M$). En effet, si $y_1$, $y_2$ sont dans $f^{-n}(x)$ et $l$ désigne le plus petit entier tel que $f^l(y_1)=f^l(y_2)$, on a $d(f^{l-1}(y_1),f^{l-1}(y_2)) \geq \delta$. Autrement dit, 
$$h_{\mbox{top}}(f) \geq \log N=\log | \mbox{deg}(f)|.$$
Quand $f$ n'est pas un revêtement, un bon choix de $x$ permet essentiellement de se ramener au cas précédent. En effet, le théorème de Sard permet de prendre $x$ de sorte que la $n$-orbite des points de $f^{-n}(x)$ transite peu dans un petit voisinage de l'ensemble critique.\\\\

Voici le plan de la démonstration de la proposition \ref{prop1}. Ici $U_{\gamma}$ désigne toujours un $\gamma$-voisinage du support de $\mu$.\\
Dans un premier temps, on va majorer le cardinal d'un ensemble $(n, \delta)$-séparé hors de $U_{\gamma}$ par $d^{n(1+\epsilon)}$. Il s'agit donc de localiser au complémentaire de $U_{\gamma}$ l'argument de Gromov. Autrement dit, on sera essentiellement ramené à majorer $\mbox{vol}(\Gamma_n | U_{\gamma}^c) = \int_{\Gamma_n | U_{\gamma}^c} \omega_n \wedge \omega_n$ par $d^{n(1+\epsilon)}$. Ce raisonnement est valable pour tout endomorphisme de $\Pp^2(\Cc)$.\\
La seconde étape consiste à construire un ensemble $(n, \delta)$-séparé contenu dans un ensemble $P_n$ de points de $f^{-n}(x)$, de cardinal minoré par $\mbox{Cardinal}(P_n) d^{-2 \alpha n}$. C'est dans cette étape que l'on utilise la généricité de $f$.

\subsection{{\bf Majoration du cardinal d'ensembles $(n, \delta)$-séparés}}

Dans ce paragraphe, on va démontrer le lemme suivant:
\begin{lemme}{\label{majoration}}
Pour $f$ un endomorphisme holomorphe quelconque, le cardinal d'un ensemble $(n, \delta)$-séparé hors de $U_{\gamma}$ (avec $\delta$ petit) est majoré par $Cnd^n$ (où $C$ est une constante qui ne dépend que de $\delta$ et $\gamma$).

\end{lemme}

 \begin{proof}
Dans la démonstration, on notera toujours $C$ toute constante qui ne dépend que de $\delta$ et $\gamma$.\\
Soit $F$ un ensemble $(n, \delta)$-séparé hors de $U_{\gamma}$. Il induit un ensemble $G$ $\delta$-séparé dans $\Gamma_n |(U_{\gamma})^c$.\\
On a alors,
$$\mbox{vol}(\Gamma_n | U_{\frac{\gamma}{2}}^c  ) \geq \sum_{y \in G} \mbox{vol}(B_n(y, \frac{\delta}{2}) \cap \Gamma_n) \geq C \mbox{ Card}(F),$$
par le théorème de Lelong.\\
On obtient donc une minoration du volume du multigraphe restreint à $U_{\frac{\gamma}{2}}^c$ par $C \mbox{ Card}(F)$ (avec $\delta$ petit).\\
Pour dominer le cardinal de $F$ par $Cnd^n$, il reste donc à majorer $\mbox{vol}(\Gamma_n | U_{\frac{\gamma}{2}}^c)$ par cette même quantité.\\
Tout d'abord, on a:
$$\mbox{vol}(\Gamma_n | U_{\frac{\gamma}{2}}^c) = \int_{\Gamma_n | U_{\frac{\gamma}{2}}^c} \omega_n \wedge \omega_n = \int_{U_{\frac{\gamma}{2}}^c} \sum_{i,j=0}^{n-1} f^{i *} \omega \wedge f^{j *} \omega.$$
La majoration de $\mbox{vol}(\Gamma_n | U_{\frac{\gamma}{2}}^c)$ par $Cnd^n$ se déduit donc de celle de 
$$\int_{U_{\frac{\gamma}{2}}^c}  \frac{f^{i *} \omega}{d^i} \wedge \frac{f^{j *} \omega}{d^j}$$
par $C(\frac{1}{d^i} +\frac{1}{d^j})$.\\
On notera $T_i=\frac{f^{i *} \omega}{d^i}$. Rappelons la construction de $T$ (voir \cite{FS}).\\
Comme la forme $f^{*}\omega$ est cohomologue à $d \omega$, on a:
$$\frac{f^{*}\omega}{d}= \omega +dd^c u,$$
où $u$ est une fonction lisse de $\Pp^2(\Cc)$.\\
En itérant cette relation, on obtient que:
$$T_i= \omega +dd^c G_i,$$
avec $G_i=\sum_{l=0}^{i-1} \frac{u \circ f^l}{d^l}$.\\
Autrement dit, en passant à la limite, on a $T= \omega + dd^c G$ où $G$ est une fontion continue qui vérifie:
$$\max_{\Pp^2(\Cc)}|G_{i}-G| \leq \frac{C}{d^i}.$$
Soit maintenant $\psi$ une fonction $C^{\infty}$ à support compact dans $U_{\frac{\gamma}{3}}^c$, comprise entre $0$ et $1$ et qui vaut $1$ sur $U_{\frac{\gamma}{2}}^c$. Alors:
$$\int_{ U_{\frac{\gamma}{2}}^c} T_i \wedge T_j=\int_{ U_{\frac{\gamma}{2}}^c} \psi T_i \wedge T_j \leq \int \psi T_i \wedge T_j.$$
Mais la dernière intégrale est égale à:
$$\int \psi (T_i-T) \wedge T_j + \int \psi T \wedge (T_j-T)$$
car $\int \psi T \wedge T$ est dominée par $T \wedge T (U_{\frac{\gamma}{3}}^c)$ qui vaut $0$.\\
Autrement dit, si on sait majorer $\int \psi (T_i-T) \wedge S$ (où $S$ désigne un $(1,1)$-courant positif fermé de masse $1$) par $\frac{C}{d^i}$, on pourra conclure.\\
Cependant, on a:
$$\int \psi (T_i-T) \wedge S=\int (G_i - G) dd^c \psi \wedge S,$$
c'est-à-dire,
$$\int \psi (T_i-T) \wedge S \leq |G_i - G|_{\Pp^2(\Cc)} |\psi|_{C^2},$$
qui est bien majoré par $\frac{C}{d^i}$.

\end{proof}

Remarquons au passage le:
\begin{corollaire*}
L'entropie topologique $h_{\mathrm{top}}(f | U_{\gamma}^c)$ est majorée par $\log d$.

\end{corollaire*}

Ici, $h_{\mbox{top}}(f | U_{\gamma}^c)$ désigne l'entropie topologique de $f$ localisée à $U_{\gamma}^c$. Elle est définie par: 
$$h_{\mbox{top}}(f|U_{\gamma}^c)=\sup_{\delta > 0} \limsup_{n \rightarrow \infty} \frac{1}{n} \log (\max \{\mbox{Card}(F),\mbox{ } F \mbox{ } (n,\delta) \mbox{ séparé}, F \subset U_{\gamma}^c \}).$$

\subsection{{\bf Construction d'ensembles $(n, \delta)$-séparés}}

On considère un $\epsilon$-voisinage, $C_{\epsilon}$, de l'ensemble critique de $f$. L'application $f$ est localement injective sur le complémentaire de $C_{\epsilon}$. Il existe donc $\delta > 0$ tel que pour $x$ et $y$ dans $(C_{\epsilon})^c$ avec $d(x,y) \leq \delta$ on ait $f(x) \neq f(y)$. Enfin, on prend $\delta$ strictement plus petit que $\epsilon$.\\
L'objectif de ce paragraphe est alors de démontrer:

\begin{lemme}
On peut construire un ensemble $(n, \delta)$-séparé contenu dans un
ensemble $P_n$ de points de $f^{-n}(x)$, de cardinal minoré par
$\mathrm{Card }(P_n) d^{-2m}$, où $m$ est le nombre maximal de passages de l'orbite d'ordre $n$ d'un point de $P_n$ dans $C_{2\epsilon}$.

\end{lemme}

\begin{proof}

La démonstration de ce lemme va se faire en deux étapes. Dans la première, on va construire un ensemble $(n, \delta)$-séparé contenu dans l'ensemble $P_n$. Dans la seconde, on minorera le cardinal de cet ensemble.\\\\
On note $P_{n-k}$  l'ensemble $f^k(P_n)$.\\
Dans un premier temps, on part de $x$ et on considère les points de $P_1 \subset f^{-1}(x)$.\\
Parmi eux, il y a ceux qui sont dans $(C_{2\epsilon})^c$ et ceux qui n'y sont pas. On va garder tous les points qui sont dans la première catégorie. Ensuite, parmi ceux qui sont dans $C_{2\epsilon}$, on ne garde que celui qui a le plus d'antécédents dans $P_n$. On note $Q_1$ l'ensemble constitué par ce point et par les points de $P_1 \cap (C_{2\epsilon})^c$. D'autre part, on dira que l'on a eu une mauvaise transition au point de $(C_{2\epsilon})$ que l'on a conservé.\\
En remplaçant maintenant $x$ par les points de $Q_1$ dans le raisonnement précédent, on construit un ensemble $Q_2 \subset P_2$. Puis, en itérant le procédé, on arrive à un ensemble $Q_n \subset P_n \subset f^{-n}(x)$.\\
L'ensemble $Q_n$ obtenu est $(n, \delta)$-séparé.\\
En effet, si $y_1,y_2 \in Q_n$ et $d(f^k(y_1),f^k(y_2)) \leq \delta$ pour $k=0,...,n-1$ alors $f^{n-1}(y_1)=f^{n-1}(y_2)$ car de deux choses l'une:\\
soit $f^{n-1}(y_1)$ et $f^{n-1}(y_2)$ sont dans $(C_{\epsilon})^c$ et l'égalité provient de la définition de $\delta$;\\
soit $f^{n-1}(y_1)$ ou $f^{n-1}(y_2)$ est dans $C_{\epsilon}$ et alors ces deux éléments sont dans $C_{2\epsilon}$ (car $\delta$ est inférieur à $\epsilon$) et l'égalité découle de la définition de $Q_1$.\\
A partir de là, on montre de même que $f^{n-2}(y_1)=f^{n-2}(y_2)$ et ainsi de suite jusqu'à $y_1=y_2$.\\
On a donc construit un ensemble $(n, \delta)$-séparé, $Q_n$, inclus dans $P_n$. Il reste à minorer le cardinal de cet ensemble par $\mbox{Card}(P_n) d^{-2m}$.\\\\
On va faire une récurrence sur le nombre maximal $k$ de mauvaises transitions entre $x$ et un point $y$ de $P_n$ (ce nombre varie entre $0$ et $m$). Elle va montrer que le cardinal de $Q_n$ est minoré par $\mbox{Card}(P_n) d^{-2k}$.\\
Pour $k=0$, le résultat est clair car on n'enlève aucun point de $P_n$.\\
On va traiter le cas $k=1$ pour mieux comprendre le procédé.\\
Si on part d'un point $y$ dans $P_n$, on a deux possibilités:\\
- Soit il existe $l$ dans $ \{1,...,n\}$ tel que $f^l (y)$ soit une mauvaise transition. Dans ce cas, on note $S(y)$ l'ensemble des préimages dans $P_n$ des points de $f^{-1}(f^{(l+1)}(y)) \cap P_{n-l} \cap C_{2\epsilon}$. Comme on ne garde que la branche qui donne le plus de points dans $P_n$, on en déduit que parmi les points de $S(y)$, il y en a au moins $\mbox{Card}(S(y))d^{-2}$ qui sont dans $Q_n$.\\
- Soit ce $l$ n'existe pas et alors $y$ est dans $Q_n$. Dans ce cas, on note $S(y)=y$.\\
Maintenant, on peut recommencer avec $z \in P_n - S(y)$ et ainsi de suite. On obtient donc:
$$\mbox{Card}(Q_n) \geq \mbox{Card}(P_n)d^{-2}.$$
On suppose la propriété vraie jusqu'au rang $k$ et on veut la montrer au rang $k+1$.\\
Quand on enlève les points de $P_n$ correspondant aux $k$ premières mauvaises transitions, on obtient un sous-ensemble de $P_n$ (qui contient $Q_n$) de cardinal minoré par $\mbox{Card}(P_n) d^{-2k}$.\\
Dans chaque branche de ce sous-ensemble il reste au plus une mauvaise transition. On se retrouve donc dans le cas  $k=1$ avec ce sous-ensemble à la place de $P_n$. Autrement dit, on trouve bien une minoration de $Q_n$ par $d^{-2} ( \mbox{Card}(P_n)d^{-2k})=\mbox{Card}(P_n)d^{-2(k+1)} $. La récurrence est donc démontrée.\\
Dans notre situation, on sait que le nombre de mauvaises transitions est majoré par $m$.\\
On obtient donc une minoration du cardinal de l'ensemble $Q_n$, $(n, \delta)$-séparé, par $\mbox{Card}(P_n) d^{-2m}$.

\end{proof}

\subsection{{\bf Fin de la démonstration de la proposition \ref{prop1} }}

On fixe $\alpha > 0$ arbitrairement petit et dans toute la suite $n$ sera supposé grand.\\
On rappelle que l'on considère des endomorphismes holomorphes génériques. Cela signifie que $C_f \cap f^i(C_f) \cap f^j(C_f)= \emptyset$ dès que $i,j \in \Nn^{*}$ et $i \neq j$. La généricité de cette condition est démontrée dans \cite{FS1}.\\
Soit $k$ tel que $\frac{3}{k} < \alpha$.\\
Tout d'abord, si $\epsilon$ est assez petit, on a $C_{2\epsilon} \cap f^i(C_{2\epsilon}) \cap f^j(C_{2\epsilon})= \emptyset$ si $i,j \in \{1,...,k\}$ et $i \neq j$.\\
Ensuite, l'application $f$ est localement injective sur le complémentaire de $C_{\epsilon}$. Il existe donc $\delta=\delta(\alpha,f) > 0$ tel que pour $x$ et $y$ dans $(C_{\epsilon})^c$ avec $d(x,y) \leq \delta$ on ait $f(x) \neq f(y)$. Enfin, on prend $\delta$ petit devant $\epsilon$ et $\gamma$.\\\\
Maintenant, en utilisant le paragraphe précédent, si on fixe un point $x$ de  $\Pp^2(\Cc)$ et que l'on note $P_n$ l'ensemble de ses préimages par $f^n$ qui se trouvent hors de $U_{\gamma}$, on peut construire un ensemble $(n, \delta)$-séparé inclus dans $P_n$ de cardinal minoré par $\mbox{Card}(P_n)d^{-2m}$ (où $m$ est le nombre maximal de passages de l'orbite d'ordre $n$ d'un point de $P_n$ dans $C_{2\epsilon}$).\\
D'autre part, en utilisant la majoration obtenue dans le lemme \ref{majoration}, le cardinal de cet ensemble est majoré par $Cnd^n$ (où $C$ est une constante qui ne dépend que de $\delta$ et $\gamma$).\\
Autrement dit, on a:
$$\mbox{Card}(P_n) \leq Cnd^{2m}d^{n}.$$\\
Maintenant, si on montre que pour tout point $y$ de $\Pp^2(\Cc)$, le nombre de passages de l'orbite d'ordre $n$ de $y$ dans $C_{2\epsilon}$ est majoré par $\alpha n$ (pour $n$ grand), on aura:
$$\mbox{Card}(P_n) \leq Cnd^{n(2\alpha + 1)},$$
et la proposition sera démontrée.\\\\
Soit donc $A$ l'ensemble des points de $\Pp^2(\Cc)$ dont l'orbite visite souvent un $2\epsilon$-voisinage de l'ensemble critique, i.e.:
$$A=\{ x \in \Pp^2(\Cc) / \mbox{ Card} (C_{2\epsilon} \cap \{x,f(x),...,f^{n-1}(x) \}) \geq \alpha n \}.$$
Alors $A$ est bien vide:\\
en effet, soit $y \in A$. Si on prend $m \in \Nn$ avec $m+k \leq n$, on a:
$$\mbox{Card}(i \in \{m,...,m+k\} \mbox{, } f^i(y) \in C_{2\epsilon}) \leq 2$$
par définition de $\epsilon$.
D'où
$$\mbox{Card}(i \in \{1,...,n\} \mbox{, } f^i(y) \in C_{2\epsilon}) \leq 2([\frac{n}{k}]+1) \leq n(\frac{2}{k}+\frac{2}{n}) < \alpha n.$$

\section{{\bf Contrôle du genre}}

Dans ce paragraphe, on va passer à la démonstration du théorème \ref{genre}. Dans celle-ci, il suffit de traiter le cas où $L$ est générique. En effet, une majoration du genre de $f^{-n}(L)-U$ pour des droites génériques par $Cd^{(1+ \alpha)n}$ (avec $C$ indépendante de $L$) conduit à la même majoration du genre de $f^{-n}(L)-U$ pour toutes les droites par passage à la limite. Dans toute la suite on supposera donc que $f^{-n}(L)$ est lisse et que $L$ est transverse à l'ensemble postcritique.\\
Pour majorer le genre de $f^{-n}(L)-U$, on va construire une partition en disques et anneaux de $f^{-n}(L)$ et faire un calcul de caractéristique d'Euler. Cependant, lors de la démonstration du cas critiquement fini, on a vu que l'on devait contrôler la longueur des arêtes de la triangulation. Dans ce but, on va utiliser ici un peu de géométrie hyperbolique. En effet, si on considère une petite arête dans la partie épaisse d'une surface hyperbolique, on peut l'insérer dans un disque de sorte à contrôler le module de l'anneau ainsi créé. On est donc en mesure d'utiliser un argument longueur-aire.\\
Voici le plan de ce paragraphe: la décomposition de $f^{-n}(L)$ en disques et anneaux occupera les deux premières parties. En effet, dans un premier temps on construira de telles partitions sur des surfaces hyperboliques quelconques d'aire finie. Puis, on utilisera cette construction de façon itérative pour produire la partition de $f^{-n}(L)$ que l'on cherche. Enfin les deux dernières parties seront consacrées d'une part à la majoration du genre et d'autre part au contrôle des longueurs des arêtes de la triangulation.

\subsection{{\bf Construction d'une partition statique}}

Dans cette partie, on va produire une partition en disques et anneaux d'une surface de Riemann hyperbolique $S$ quelconque d'aire finie.\\
Rappelons que $S$ est composée d'une partie mince (ensemble des points où le rayon d'injectivité est inférieur à $\mbox{Argsh}(1)$) et d'une partie épaisse. Par ailleurs la partie mince consiste en cusps et anneaux (voir \cite{Bu} pour plus de détails).\\
La partition de $S$ sera donc composée de la partie mince et d'un recouvrement régulier de la partie épaisse par des disques.\\
On note $F$ un ensemble $\epsilon$-séparé maximal dans la partie épaisse (avec $\epsilon < \mbox{Argsh}(1)$). L'aire pour une métrique de courbure $-1$ d'un disque $D(x,\frac{\epsilon}{2})$ avec $x$ dans la partie épaisse est une constante qui ne dépend que de $\epsilon$. Le cardinal de $F$ est donc majoré par $C \mbox{aire}(S)$ (on note $C$ toute constante qui ne dépend que de $\epsilon$). On recouvre la partie épaisse avec la réunion des disques $D(x, \epsilon)$ (où $x$ décrit $F$). On en déduit une partition recouvrant la partie épaisse par des disques qui sont les intersections des disques précédents, d'où une partition de $S$ par des disques et des anneaux.\\
Si $x$ est un élément de $F$, il y a au plus $C$ points de $F$ dans
$D(x,2 \epsilon)$. La partition ci-dessus contient donc au plus $C
\mbox{aire}(S)$ arêtes. Ce nombre prend en compte les arêtes qui bordent la partie mince.\\
C'est cette partition que l'on utilisera pour construire une partition dynamique de $f^{-n}(L)$.\\
Remarquons enfin que l'aire de $S$ est bornée par sa topologie (grâce à la formule de Gauss-Bonnet). Autrement dit, le nombre d'arêtes de la partition précédente est majoré par $-C \chi(S)$.

\subsection{{\bf Construction de la partition dynamique de $f^{-n}(L)$}}

On part d'une droite $L$ projective et on note $V$ l'ensemble des valeurs critiques de $f$.\\
La première étape consiste à fabriquer une ``bonne'' partition de $L$ composée de disques et d'anneaux.\\
Quitte à bouger un peu $L$, $S=L - L \cap V$ est une surface de Riemann hyperbolique (car $V$ est de degré $3d(d-1) \geq 3$). En utilisant le paragraphe précédent, on peut donc construire une partition de $S$, donc de $L$, en des disques et des anneaux. Elle contient au plus $C \mbox{Card}( L \cap V)$ arêtes.\\
C'est la ``bonne'' partition que l'on cherchait.\\
Par construction, les disques (resp. les anneaux) de la partition de $L$ se relèvent sous forme de disques (resp. d'anneaux). En effet, l'application $f$ est un revêtement de $f^{-1}(L) - f^{-1}(V)$ sur $L - V$. Ainsi, en tirant en arrière par $f$ la partition de $L$, on en obtient une de $f^{-1}(L)$.\\
On va maintenant raisonner dans $f^{-1}(L)$ où on veut construire une ``bonne'' partition. Quitte à bouger un peu les partitions, on supposera dans la suite que les arêtes que l'on construit ne touchent jamais un itéré de $V$.\\
Les disques (resp. les anneaux) de $f^{-1}(L)$ qui rencontrent $V$ en au plus un point (resp. qui ne rencontrent pas $V$) se relèvent par $f$  sous forme de disques (resp. d'anneaux). On ne va donc pas les modifier. Par contre, on va redécouper les disques et anneaux qui n'entrent pas dans les catégories ci-dessus.\\
On va traiter le cas du disque (celui de l'anneau est identique).\\
On note $D$ un disque de $f^{-1}(L)$ qui touche $V$ en au moins deux points. En doublant $D$ et en enlevant les points de $V \cap D$ ainsi que leurs symétriques, on obtient une surface de Riemann $S^{'}$ hyperbolique.\\
Comme dans l'étape précédente, on produit une partition de $S^{'}$ en utilisant les anneaux de la partie mince auxquels on ajoute un recouvrement de la partie épaisse par des disques $D(x, \epsilon)$.\\
Pour des raisons de symétrie, le bord de $D$ dans $S^{'}$ est une géodésique. La partition précédente induit donc une partition de $D - D \cap V$, donc de $D$, en disques et anneaux qui contient $C \mbox{Card}( D \cap V)$ arêtes.\\
En recommençant ce que l'on vient de faire avec les autres disques et anneaux qui rencontrent $V$, on construit une ``bonne'' partition de $f^{-1}(L)$. Si on relève celle-ci par $f$ et que l'on itère le procédé, on aboutit à une partition ``dynamique'' de  $f^{-n}(L)$ en disques et anneaux.\\
Dans le paragraphe suivant, on va montrer qu'un contrôle des longueurs
des arêtes de la partition induit la majoration du genre de
$f^{-n}(L)$ hors de $U_{\gamma}$ en $d^{n(1+\alpha)}$.

\subsection{{\bf Majoration du genre des préimages de droites}}

On va procéder ici comme dans le cas critiquement fini.\\
Soit $\Gamma$ la réunion des arêtes et des sommets de la partition ainsi obtenue qui sont contenus dans $U_{\gamma}$.\\
Alors:
$$g(f^{-n}(L) - U_{\gamma}) \leq g(f^{-n}(L) - \Gamma),$$
et c'est ce dernier terme que l'on va majorer.\\
Une composante connexe de $f^{-n}(L) - \Gamma$ peut être de deux types. Soit elle ne contient pas d'arête de la triangulation et alors c'est nécessairement un disque ou un anneau (donc de genre nul). Soit elle en contient au moins une. Si $\mathcal{C}$ désigne l'union des composantes du deuxième type, on obtient :
$$g(f^{-n}(L) - U_{\gamma}) \leq \mbox{Nombre de composantes de } \mathcal{C} - \frac{ \chi(\mathcal{C})}{2},$$
Ainsi:
$$ g(f^{-n}(L) - U_{\gamma}) \leq \frac{3}{2} \mbox{du nombre d'arêtes
  qui sortent de }U_{\gamma}.$$
En effet les faces qui composent les éléments de $\mathcal{C}$ étant des disques ou des anneaux, $- \chi(\mathcal{C})$ est majoré par le nombre d'arêtes qui sortent de $U_{\gamma}$.\\
Il reste donc à majorer ce dernier terme.\\
Si $a$ est une arête qui sort de $U_{\gamma}$, elle entre dans un des deux cas suivants:\\
\underline{1\up{er} cas}: $a$ possède un sommet hors de $U_{\frac{\gamma}{2}}$.\\
La majoration du nombre de ces arêtes découle alors du contrôle des préimages des sommets.\\
On fixe $\alpha > 0$.\\
Il existe $n_0$ à partir duquel on a:
$$\max_{x \in \Pp^2} \mbox{Card} (f^{-n}(x) \cap U_{\frac{\gamma}{2}}^c) \leq d^{n(1+\alpha)}.$$
Etant donné que le nombre d'arêtes créées au rang $k$ est de l'ordre de $Cd^k$, on obtient alors une majoration du nombre d'arêtes de $f^{-n}(L)$ qui ont un sommet hors de $U_{\frac{\gamma}{2}}$ par:
$$ \sum_{k=0}^{n-n_0} Cd^k d^{(n-k)(1+ \alpha)} + \sum_{k=n-n_0+1}^{n} Cd^k d^{2(n-k)}.$$
Le nombre d'arêtes de la partition de $f^{-n}(L)$ qui entrent dans le premier cas est donc majoré par $d^{n(1+2 \alpha)}$ (pour $n$ grand).\\
\underline{2\up{ème} Cas}: Les deux sommets de $a$ sont dans $U_{\frac{\gamma}{2}}$ (en particulier la longueur de $a$ est supérieure à $\gamma$).\\
Une majoration du nombre des arêtes de ce type par $d^{n(1+\alpha)}$ démontrerait le résultat.\\
C'est ce contrôle qui va être l'objet du paragraphe suivant.

\subsection{{\bf Contrôle des longueurs des arêtes de la partition}}

On va voir que le contrôle des longueurs des arêtes de la partition de $f^{-n}(L)$ est possible quitte à la bouger un peu.
Soit $D$ un disque de $f^{-k}(L)$ qui touche $V$ en au moins deux points (le cas de l'anneau est identique). Rappelons que l'on note $S^{'}$ la surface de Riemann hyperbolique associée à $D$ et $\mu$ sa métrique de Poincaré.\\
Parmi les éléments qui partitionnent $D$, on a des disques $D(x,\epsilon) \cap D$.\\
On va voir que quitte à bouger un peu $\partial D(x,\epsilon)$, on peut contrôler les longueurs des préimages par $f^{n-k}$ de $\partial D(x,\epsilon) - \{1 \mbox{ point} \}$.\\
Traitons l'exemple d'un disque $D(x,\epsilon) \cap D$ où $x$ est à distance au moins $2 \epsilon$ de $\partial D$.\\
Sur le disque $\Delta$, constitué de l'anneau $D(x,2 \epsilon) - D(x,\epsilon)$ auquel on a enlevé une géodésique $c$ qui joint $\partial D(x,2 \epsilon)$ à $\partial D(x,\epsilon)$, on peut définir $d^{2(n-k)}$ branches inverses $f_j^{-n+k}$ de $f^{n-k}$.\\
D'autre part, si on note $ds^2$ la métrique ambiante restreinte à $f^{-n}(L)$, $(f_j^{-n+k})^{*} ds^2$ est conforme à la métrique hyperbolique (écrite en coordonnées polaires) $\mu=dr^2 + \mbox{sh}^2(r) d \theta^2$ sur $\Delta$. Cela se traduit par:
$$(f_j^{-n+k})^{*} ds^2= \rho_j^2 (dr^2 + \mbox{sh}^2(r) d \theta^2).$$
Voici l'argument longueur-aire donnant le contrôle voulu:
$$\mbox{aire}(f^{-n+k}(\Delta))=\int_0^{2 \pi} \int_{\epsilon}^{2\epsilon} \sum_{j=1}^{d^{2(n-k)}} \rho_j^2  
\mbox{sh}(r) dr d\theta \leq \mbox{aire}(f^{-n+k}D(x,2 \epsilon))=A.$$
Il existe alors $\epsilon_0$ dans l'intervalle $[\epsilon, 2 \epsilon]$ qui vérifie:
$$\sum_{j=1}^{d^{2(n-k)}} \int_0^{2 \pi} \rho_j^2 (\epsilon_0,.) d\theta \leq C A.$$
On a donc:
$$\mbox{Long}(f_j^{-n+k}( \partial D(x,\epsilon_0) - c))= \int_0^{2 \pi} \rho_j(\epsilon_0,.) \mbox{sh}(\epsilon_0) d\theta \leq C\left( \int_0^{2 \pi} \rho_j^2 (\epsilon_0,.) d\theta \right)^{\frac{1}{2}} \leq \gamma$$
sauf un nombre d'indices $j$ majoré par $CA$.\\
Pour la partition de $D$, on remplace alors $D(x,\epsilon)$ par $D(x,\epsilon_0)$.\\
Les disques de taille $2\epsilon$ se recouvrent un nombre fini borné à priori de fois. Autrement dit, les arêtes créées au rang $k$, mais non incluses dans les préimages des arêtes des étapes précédentes, ont leurs préimages dans $f^{-n}(L)$ de longueur inférieure à $\gamma$ (sauf $Cd^n$ d'entre elles).
Ainsi, en considérant maintenant toutes les étapes $k$ et en remarquant que redécouper une arête créée au rang $k$ dans une étape postérieure ne change rien aux estimées, on en déduit que le nombre d'arêtes de $f^{-n}(L)$ de longueur supérieure à $\gamma$ est dominé par $Cnd^n$.

\bigskip

Henry de Thélin\\
Université Paris-Sud (Paris 11), Mathématique, Bât. 425,
91405 Orsay, France.\\
E-mail: Henry.De-Thelin@math.u-psud.fr

\end{document}